\documentclass[11pt]{article}
 \usepackage{amsfonts,amssymb,enumerate,epsf,dsfont,amsmath,array}
 \usepackage{graphicx}
\date{}
\usepackage{hyperref}
\usepackage{xcolor}
\usepackage{authblk}
\usepackage{standalone}

\title{Asymptotic normality of degree counts}
\usepackage{verbatim}
 
 \usepackage[english]{babel}
 \def\botcaption#1#2{\medskip\centerline{{\scshape #1.}\kern8pt
 {\rm #2}}\bigskip}

 \makeatletter
 \@addtoreset{equation}{section}
 \makeatother

 \newtheorem{ittheorem}{Theorem}
 \newtheorem{itlemma}{Lemma}
 \newtheorem{itproposition}{Proposition}
 \newtheorem{itdefinition}{Definition}
 \newtheorem{itremark}{Remark}
 \newtheorem{itclaim}{Claim}
 \newtheorem{itcorollary}{Corollary}

 \newenvironment{theorem}{\addtocounter{equation}{1}
 \begin{ittheorem}}{\end{ittheorem}}

 \newenvironment{lemma}{\addtocounter{equation}{1}
 \begin{itlemma}}{\end{itlemma}}

 \newenvironment{proposition}{\addtocounter{equation}{1}
 \begin{itproposition}}{\end{itproposition}}

 \newenvironment{definition}{\addtocounter{equation}{1}
 \begin{itdefinition}}{\end{itdefinition}}

 \newenvironment{remark}{\addtocounter{equation}{1}
 \begin{itremark}}{\end{itremark}}

 \newenvironment{claim}{\addtocounter{equation}{1}
 \begin{itclaim}}{\end{itclaim}}

\newenvironment{corollary}{\addtocounter{equation}{1}
	\begin{itcorollary}}{\end{itcorollary}}

 \newenvironment{proof}{\noindent {\bf Proof.\,}
 }{\hspace*{\fill}$\qed$\medskip}

 \newenvironment{proof*}{\noindent {\bf Proof\,}
}{\hspace*{\fill}$\qed$\medskip}

 \newcommand{\be}[1]{\begin{equation}\label{#1}}
 \newcommand{\ee}{\end{equation}}

 \newcommand{\bl}[1]{\begin{lemma}\label{#1}}
 \newcommand{\el}{\end{lemma}}

 \newcommand{\br}[1]{\begin{remark}\label{#1}}
 \newcommand{\er}{\end{remark}}

 \newcommand{\bt}[1]{\begin{theorem}\label{#1}}
 \newcommand{\et}{\end{theorem}}

 \newcommand{\bd}[1]{\begin{definition}\label{#1}}
 \newcommand{\ed}{\end{definition}}

 \newcommand{\bcl}[1]{\begin{claim}\label{#1}}
 \newcommand{\ecl}{\end{claim}}

 \newcommand{\bp}[1]{\begin{proposition}\label{#1}}
 \newcommand{\ep}{\end{proposition}}

 \newcommand{\bc}[1]{\begin{corollary}\label{#1}}
 \newcommand{\ec}{\end{corollary}}

 \newcommand{\bpr}{\begin{proof}}
 \newcommand{\epr}{\end{proof}}

 \newcommand{\bi}{\begin{itemize}}
 \newcommand{\ei}{\end{itemize}}

 \newcommand{\ben}{\begin{enumerate}}
 \newcommand{\een}{\end{enumerate}}

\newenvironment{sistema}% 
{\left\lbrace\begin{array}{@{}l@{}}}% 
	{\end{array}\right.}

 \def \ba {\begin{array}}
 \def \ea {\end{array}}

 \def \qed {{\square\hfill}}
 
 \def \R {{\mathbb R}}
 
 \def \N {{\mathbb N}}
 \def \P {{\mathbb P}}
 \def \E {{\mathbb E}}
 \def \G {{\Gamma}}
 \def \i {{\mathds 1}}
 \def \cF {{\cal F}}
 \def \d {{\delta}}

 \def \ra {\rightarrow}

\begin{document}

 	\title{Asymptotic normality of degree counts in a general preferential attachment model}
 	
 \author[
 {}\hspace{0.5pt}\protect\hyperlink{hyp:email1}{1},\protect\hyperlink{hyp:affil1}{a}%,\protect\hyperlink{hyp:affil1}{b}
 ]
 {\protect\hypertarget{hyp:author1}{Simone Baldassarri}}
 
 \author[
 {}\hspace{0.5pt}\protect\hyperlink{hyp:email2}{2},\protect\hyperlink{hyp:affil1}{a},\protect\hyperlink{hyp:corresponding}{$\dagger$}
 ]
 {\protect\hypertarget{hyp:author2}{Gianmarco Bet}}

 \affil[ ]{\centering
 	\small\parbox{365pt}{\centering
 		\parbox{5pt}{\textsuperscript{\protect\hypertarget{hyp:affil2}{a}}}Dipartimento di Matematica e Informatica ``Ulisse Dini", Università degli Studi di Firenze, Firenze, Italy
 	}
 }
 
 \affil[ ]{\centering
 	\small\parbox{365pt}{\centering
 		\parbox{5pt}{\textsuperscript{\protect\hypertarget{hyp:email1}{1}}}\texttt{\footnotesize\href{mailto:simone.baldassarri@unifi.it}{simone.baldassarri@unifi.it}},
 		\parbox{5pt}{\textsuperscript{\protect\hypertarget{hyp:email2}{2}}}\texttt{\footnotesize\href{mailto:gianmarco.bet@unifi.it}{gianmarco.bet@unifi.it}},
 	}
 }
 
 \affil[ ]{\centering
 	\small\parbox{365pt}{\centering
 		\parbox{5pt}{\textsuperscript{\protect\hypertarget{hyp:corresponding}{$\dagger$}}}Corresponding author
 	}
 }

 	\maketitle
 	
 		\vspace{-1.25cm}
 	
 	\begin{center}
 		{\it Both authors are deeply grateful to their former group leader \\ Francesca Romana Nardi for her guidance and friendship, and for encouraging the collaboration that led to this manuscript. \\ This work is dedicated to her memory.}		
 	\end{center}
 	
\begin{abstract}
We consider the preferential attachment model. This is a growing random graph such that at each step a new vertex is added and forms $m$ connections. The neighbors of the new vertex are chosen at random with probability proportional to their degree. It is well known that the proportion of nodes with a given degree at step $n$ converges to a constant as $n\rightarrow\infty$. The goal of this paper is to investigate the asymptotic distribution of the fluctuations around this limiting value. We prove a central limit theorem for the joint distribution of \textit{all} degree counts. In particular, we give an explicit expression for the asymptotic covariance. This expression is rather complex, so we compute it numerically for various parameter choices. We also use numerical simulations to argue that the convergence is quite fast. The proof relies on the careful construction of an appropriate martingale.
\end{abstract}

 

\section{Introduction}
The unprecedented growth in size and complexity of social and economic networks in the last two decades has sparked considerable interest in the understanding of the fundamental properties of such networks. In this context, the preferential attachment (PA) model, introduced in \cite{BA}, is a well-known model for a network that grows in time. More precisely, the PA model consists of a sequence of graphs of increasing size such that each graph is obtained from the previous according to a certain probabilistic rule. Namely, at each step a new vertex is added and it forms connections with the vertices in the graph in such a way that connections with vertices having larger degrees are more likely. In fact, several different PA models appear in the literature, depending on the concrete details of the attachment mechanism. For instance, in \cite{{BBCR},{GVCH},{SRTDWW},{WR},{ZM}} the authors investigate a {\it directed} PA model, while in \cite{{BBCS},{CF},{GV},{GVCH},{Mori},{RS}} an {\it undirected} version is considered. However, classical PA networks do not always fit real-world network data well, or in many applications it is natural to assign some kind of features to the vertices or to the edges. This led to consider some extensions of the classical PA model. For example, in \cite{BR} the authors consider a general family of preferential attachment models with multi-type edges, while \cite{{AHRS},{M},{PSY}} investigate a PA model which mixes PA rules with uniform attachment rules.

Here, we consider the PA model without self-loops described for example in \cite{GV}. In particular, each graph in the sequence is connected. We make this choice in order to simplify the calculations, but we believe that our result also holds, for example, for the PA model with self-loops considered in the standard reference \cite[Chapter 8]{VdH}. We argue this in more detail later on. We assume that each new node has a fixed number $m\in\mathbb{N}$ of edges attached to it. This is a particular case of the model considered in \cite{GV}, where $m$ is a random variable and is sampled for each new node. Our main result is a central limit theorem for the proportion of nodes with a given degree. In fact, we prove this \textit{jointly} for all degree counts. In particular, we give an explicit expression for the asymptotic covariance. The first results concerning the study of the asymptotic normality of degree counts in the preferiantal attachment models without self-loops  is given in \cite{Mori} by using martingale central limit theorems. Our results generalize those obtained in \cite{RS} for the preferential attachment tree. More precisely, in \cite{RS} the authors consider a PA model with self-loops and such that $m=1$. Note that here we consider the PA model with $m\geq1$ and without self-loops. However, this does not influence the asymptotic behavior of the degree counts, since as the graph size goes to infinity, the probability that a new vertex forms a self-loop tends to zero. Because of this one would expect to recover the results in \cite{RS} when plugging $m=1$ in our result. Indeed this is the case if one takes into account a few minor mistakes in \cite{RS} which we will discuss later. Note that a major difference between the two models is the resulting connectivity structure. Our model produces a connected graph with probability $1$ (w.p.$1$), while the model in \cite{RS} is disconnected w.p.$1$. However, this does not play a role in the distribution of the degree counts. In \cite{PRR} the authors studied the joint degree counts in linear preferential attachment random graphs. The results are stated in terms of weight of vertices, but they can be tought of as degree, since each time a vertex receives a new edge, its weight increases by one. The main difference between our and their model is that we consider the attachment probabilities proportional to a linear function of the degree of an old vertex (see \eqref{prob}), while in \cite{PRR} they are proportional to the degree of an old vertex.

In practice, not all nodes that enter the network have the same degree, and thus it would be interesting to extend our result to the case of a random initial degree distribution. Promising results on this model have been obtained in \cite{{DEH},{GV}}. Moreover, in this paper we assume that the parameters of the model are known, but in many practical situations one is given a realization of the graph and the task is estimating the unknown parameters, see \cite{{GRW},{WWDR},{WWDR2}}. If we consider a more general class of preferential attachment graphs, for which a model-free approach is used and therefore the exact distribution of the graph is not known (see for instance \cite{KP}), we expect that the techniques presented in this paper could be used to derive central limit theorem for all the degree counts. This is an interesting open problem.

The outline of the paper is as follows: in Section \ref{res} we define the model rigorously and we state our main result. In Section \ref{proof} we give the proof of our main result.

\section{Model and results}
\label{res}
 	
Let us now describe in detail the random graph model that we consider. Fix once and for all an integer $m\geq1$. Formally, the preferential attachment model is a sequence of random graphs $(\hbox{PA}_s)_{s=1}^{t}$. The index $s$ is interpreted as a time parameter. At time $s$, the graph $\hbox{PA}_s$ has a set $V=\{0,1,...,s\}$ of $s+1$ vertices. For $s=1$ the graph PA$_1$ consists of the vertices 0 and 1, connected by $m$ edges. For $s\geq2$, the graph PA$_s$ is obtained from PA$_{s-1}$ by adding a new vertex $s$ with degree $m$ as follows. Define PA$_{s,0}=$PA$_{s-1}$ and PA$_{s,1}$,...,PA$_{s,m}$ as the intermediate graphs obtained by adding a new edge sequentially to PA$_{s,0}$. For $i=1,...,m$, PA$_{s,i}$ is obtained from PA$_{s,i-1}$ by drawing an edge from $s$ to a randomly selected vertex among $\{0,1,...,s-1\}$. The probability that a vertex $s$ is connected to some vertex $i$ is proportional to the degree of $i$. In other words, vertices with large degrees are more likely to attract new edges. We denote by $N_k(s,i)$ for $i= 1, . . . ,m$ the number of vertices of degree $k$ after the $i$-th edge has been attached at time $s$, excluding the vertex $s$. We set $N_k(s+1,0):=N_k(s,m)$. Furthermore, we denote by $D_{s,i}$ the degree of the vertex which has been attached to the $i-$th edge added when constructing PA$_s$ from PA$_{s-1}$. Consider now the $\sigma$-algebra $\cF_{s,i}$ generated by the preferential attachment construction up until the attachment of the $i$-th edge of the new vertex at time $s$. The conditional probability that the $i$-th edge connects to a vertex of degree $D_{s,i}$ is
 	
\be{prob}
\P(D_{s,i}=k|\cF_{s,i-1})=\frac{(k+\d)N_k(s,i-1)}{\sum_j (j+\d)N_j(s,i-1)},
\ee

\noindent
where $\d>-m$ is an affine parameter. The normalizing constant in (\ref{prob}) takes the simple form \cite{GV}
 	\be{}
 	\sum_{j=1}^\infty (j+\d)N_j(s,i-1)= s(2m+\d)-2m+i-1.
 	\ee
For the standard PA model considered in \cite[Chapter 8]{VdH} it is shown that there exists a probability mass function $\{p_k, \ k\geq m\}$ such that, uniformly on $i\in{\{0,...,m\}}$,
 	\be{lgn}
 	\lim_{s\ra\infty}\frac{N_k(s,i)}{s}=p_k\in{(0,1)},
 	\ee
 	almost surely, where $p_k$ is given by 
 	\be{pk}
 	p_k=(2+\d/m)\frac{\G(k+\d)\G(m+2+\d+\d/m)}{\G(m+\d)\G(k+3+\d+\d/m)}.
 	\ee	
 	Here $\G(\cdot)$ is the Gamma function. When the graph size goes to infinity, the probability that a new vertex forms a self-loop tends to zero and thus it easy to check that \eqref{lgn} and \eqref{pk} hold still for our model without self-loops following the proof proposed in \cite[Section 8.6]{VdH}. Here we choose to update the degrees during the attachment of a new vertex, but it is possible to consider also the case in which we update the degrees of the vertices only when the $m$-th edge is added. In this case, after constructing a suitable martingale with respect to the filtration $(\cF_s)_{s\geq1}$ generated by the construction of the preferential attachment graph until time $s$, it is possible to reproduce the same computations. Thus we are able to prove a similar result for this model by using the techniques presented here.

In order to state our main result, we require some further notation.
We say that the events $(A_n)_n$ hold with high probability when $P(A_n)\ra 1$ as $n\ra\infty$. Given a random vector $(X_1^{(n)}, X_2^{(n)},...)$, we write $(X_1^{(n)}, X_2^{(n)},...)\Rightarrow(X_1, X_2,...)$ to indicate that for any $k\in\N$, as $n\ra\infty$ $(X_1^{(n)}, X_2^{(n)},...,X_k^{(n)})$ converges to $(X_1, X_2,...,X_k)$ in distribution as vectors in $\R^n$.

 Our main result is the following theorem.
 	
 	\bt{nodes}
 	As $s\ra\infty$,
 	\be{gauss}
 	\Bigg(\sqrt{s}\Big(\frac{N_k(s,i)}{s}-p_k\Big), \ k=m,m+1,...\Bigg)\Rightarrow (Z_k, \ k=m,m+1,...),
 	\ee
 	where $(Z_k, \ k=m,m+1,...)$ is a mean zero Gaussian process with covariance function $R_Z$ given by
 	\be{cov}
 	\ba{lll}
 	R_Z(r,l)\\
 	=\frac{(-1)^{r+l}}{\G(l+r+1+2(\d+1)+\d/m)}\displaystyle\sum_{q=m}^{\infty}(q+\d)p_q\Big(\frac{b_m^{(l)}b_m^{(r)}}{m}\G(2(m+\d+1)+\d/m)(l+r-2m)!\\
 	\quad\quad+\dfrac{(-1)^{m+q+1}}{m}(b_m^{(l)}b_q^{(r)}+b_m^{(r)}b_q^{(l)})\G(m+q+2(\d+1)+\d/m)(l+r-q-m)!\\
 	\quad\quad+\dfrac{(-1)^{m+q+1}}{m}(b_m^{(l)}b_{q+1}^{(r)}+b_m^{(r)}b_{q+1}^{(l)})\G(m+q+1+2(\d+1)+\d/m)(l+r-m-q-1)!\\
 	\quad\quad+b_q^{(l)}b_q^{(r)}\G(2q+2(\d+1)+\d/m)(l+r-2q)!\\
 	\quad\quad+(b_q^{(l)}b_{q+1}^{(r)}+b_q^{(r)}b_{q+1}^{(l)})\G(2q+1+2(\d+1)+\d/m)(l+r-2q-1)!\\
 	\quad\quad+b_{q+1}^{(l)}b_{q+1}^{(r)}\G(2(q+\d+2)+\d/m)(l+r-2q-2)!\Big)\\
 	\quad\quad-\frac{(2m+\d)\G(l+\d)\G(r+\d)}{(\G(m+\d))^2}\displaystyle\sum_{t_1=m}^{l}\displaystyle\sum_{t_2=m}^{r}(-1)^{t_1+t_2}\\
 	\quad\quad\times\frac{[m(t_1+t_2)+2m(\d+1)+\d]^{-1}(t_1+2+\d-t_1/m)(t_2+2+\d-t_2/m)}{(t_1-m)!(t_2-m)!(l-t_1)!(r-t_2)!(t_1+2+\d+\d/m)(t_2+2+\d+\d/m)}\\
 	\quad\quad-\frac{(2m+\d)(m-1)\G(l+\d)\G(r+\d)}{m^2(\G(m+\d))^2}\displaystyle\sum_{t_1=m}^{l}\displaystyle\sum_{t_2=m}^{r}(-1)^{t_1+t_2}\\
 	\quad\quad\times\frac{(\d+t_1)(\d+t_2)[m(t_1+t_2)+2m(\d+1)+\d]^{-1}}{(t_1-m)!(t_2-m)!(l-t_1)!(r-t_2)!(t_1+2+\d+\d/m)(t_2+2+\d+\d/m)},\\
 	\ea
 	\ee
 	for $r,l\geq m$. Here $b_j^{(k)}$ are given by	
 	\be{}
 	b_j^{(k)}=\prod_{t=j}^{k-1}\frac{t+\d}{t-k}=(-1)^{k-j}\frac{\G(k+\d)}{(k-j)!\G(j+\d)}, \quad 1\leq j\leq k.
 	\ee
 	\et
 	
 	 \begin{figure}
 		\begin{center}
 			\adjustbox{max width=\textwidth}{\includestandalone[mode=image|tex]{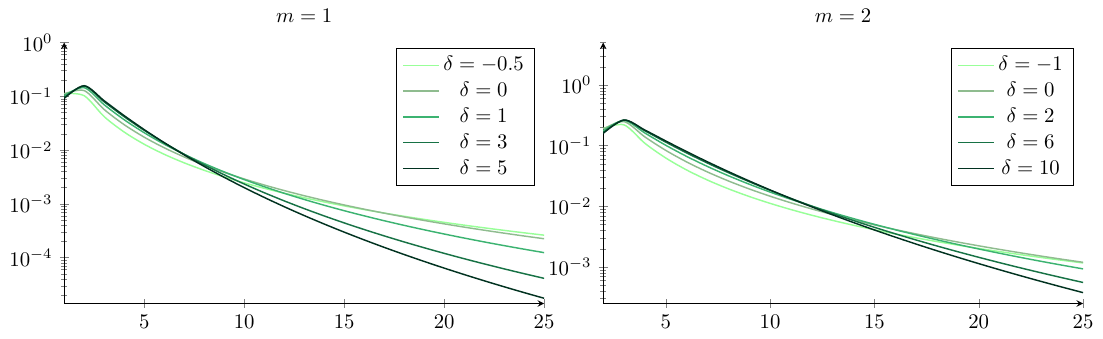}}
 		\end{center}
 		\caption{$R_Z(r,r)$ for $m=1$ and $\d=-0.5$, $0$, $1$, $3$, $5$ on the left-hand side and $R_Z(r,r)$ for $m=2$ and $\d=-1$, $0$, $2$, $6$, $10$ on the right-hand side. To highlight the different behaviour of the various variance functions, we use the logarithmic axis on the $y$-axis.}
 		\label{fig:fig1}
 	\end{figure}

 \noindent
When $m=1$, \eqref{cov} reduces to \cite[(4.28)]{RS}, but their covariance functions contains a few minor mistakes. More precisely, in the double sum the term $(2+\d)^2$ should be $(2+\d)^3$ and the terms $\G(l_1+2+2\d)$ and $\G( l_2+2+2\d)$ should be $(l_1+2+2\d)$ and $(l_2+2+2\d)$ respectively.
 
\begin{figure}[h!]
	\begin{center}
		\includestandalone[mode=image|tex]{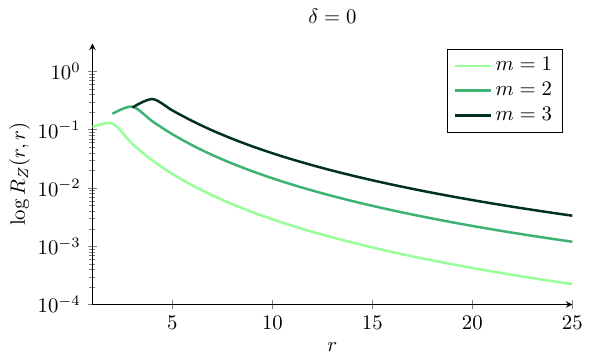}
	\end{center}
	\caption{$R_Z(r,r)$ for $\d=0$ and $m=1,2,3$. To highlight the different behaviour of the various variance functions, we use the logarithmic axis on the $y$-axis.}
	\label{fig:fig2}
\end{figure}

\begin{figure}[h!]
	\begin{center}
		\adjustbox{max width=\textwidth}{\includestandalone[mode=image|tex]{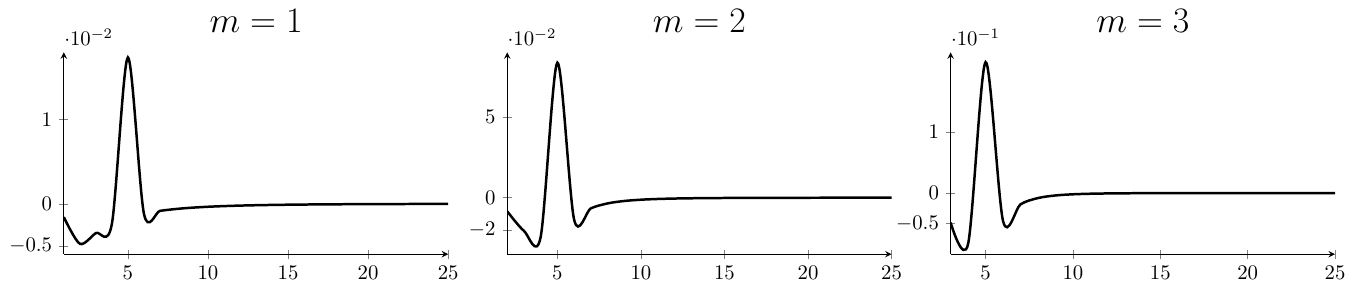}}
	\end{center}
	\caption{$R_Z(r,5)$ for $\d=0$ and $m=1$, $2$, $3$.}
	\label{fig:fig3}
\end{figure}

\begin{figure}
	\begin{center}
		\adjustbox{max width=\textwidth}{\includestandalone[mode=image|tex]{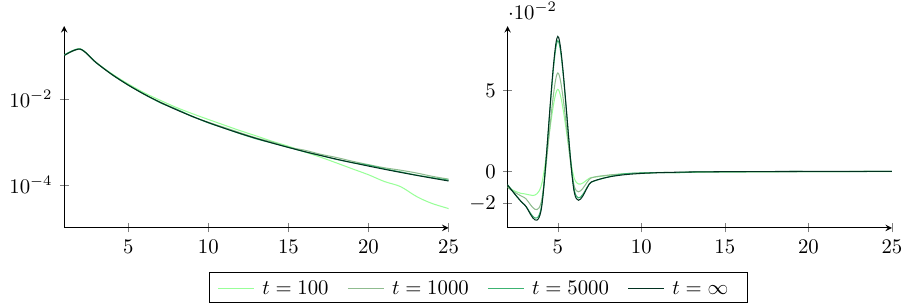}}
	\end{center}
	\caption{Left-hand side: $R_Z(r,r)$ for $\d=1$ and $m=1$ ($t=\infty$) compared with the numerical simulations stopped at time $t=100,1000,5000$. Each empirical curve was obtained by taking the average of $N=10000$ simulations. To highlight the different behaviour of the various variance functions, we use the logarithmic axis on the $y$-axis. Right-hand side: $R_Z(r,5)$ for $m=2$ and $\d=0$ compared with the numerical simulations stopped at time $t=100,1000,5000$. Each empirical curve was obtained by taking the average of $N=10000$ simulations.}
	\label{fig:fig4}
\end{figure}

Since the expression in \eqref{cov} is remarkably complicated, we plot it for various parameter choices to help the understanding. In Figure \ref{fig:fig1} we plot the function $r\mapsto R_Z(r,r)$ for fixed $m=1$ (resp.\ $m=2$) and various values of $\d$. On the other hand, in Figure \ref{fig:fig2} we plot the function $R_Z(r,r)$ for fixed $\d=0$ and various values of $m$. In Figure \ref{fig:fig3} we plot the function $R_Z(r,5)$ for $\d=0$ and $m=1$, $2$, $3$, $4$. Finally, in Figure \ref{fig:fig4} on the left-hand- side we compare the asymptotic covariance function $R_Z(r,r)$ for fixed $m=1,\d=1$ with the empirical variance obtained by numerically simulating the PA model up until time $t=100,1000,5000$. Furthermore, in Figure \ref{fig:fig4}  on the right-hand side, we compare the asymptotic covariance function $R_Z(r,5)$ for fixed $m=2,\d=0$ with the empirical covariance obtained by numerical simulation of the PA model up until time $t=100,1000,5000$. While we do not have rigorous results on the convergence speed of the rescaled vertex counts, Figure \ref{fig:fig4} suggests that the convergence is indeed quite fast.

\br{}
As explained in \cite[Chapter 8]{VdH}, it is possible to define the preferential attachment model with $m>1$ in terms of the model with $m=1$ by collapsing vertices, thus one could be tempted to directly apply this construction to the results derived in \cite{RS} for the model with $m=1$. This is a possible approach which presents its own difficulties and now we try to highlight them. The central limit theorem that we want to prove involves the number of vertices with a fixed degree, thus we need to find a relation between that quantity for the model $m>1$ and $m=1$ in order to use the result obtained in \cite{RS}. This is not straightforward, indeed a rich control over the graph construction is required at each step.
\er

\noindent
Following the argument carried out in \cite{Mori}, we are able to prove a central limit theorem for the vector composed by the rescaled number of vertices with degree greater than $k$. In this case, the covariance matrix of the limiting normal law becomes simple and now we compute it. Define the number of vertices with degree greater than $k$ as
\be{}
\psi_k(s,i):=\displaystyle\sum_{j\geq k}N_j(s,i)
\ee

\noindent
and
\be{}
\psi(s,i):=(\psi_m(s,i),\psi_{m+1}(s,i),...).
\ee

\noindent
Note that we can write
\be{}
\ba{ll}
\psi(1,i)=(0,...,0), \\
\psi(s,i+1)=\psi(s,i)+\delta_{D_{s,i}},
\ea
\ee

\noindent
where $\delta_{D_{s,i}}$ is the vector with all the coordinates equal to 0 except the ${D_{s,i}}$-th one, that is equal to 1. From \eqref{lgn} it follows that the differences of the vector valued process $\{\psi(s,i), \ s\geq1, \ i=1,...,m\}$ become more and more indipendent and identically distributed in the limit as $s\ra\infty$. Let $\pi_k(s,i):=\P(D_{s,i}=k|\cF_{s,i-1})$ and
\be{} 
\pi_k :=\lim_{s\ra\infty}\pi_k(s,i)=\dfrac{k+\d}{2m+\d}p_k.
\ee

\noindent
Thus, as $s\ra\infty$, it holds true that
\be{}
\Bigg(\sqrt{s}\Big(\frac{\psi_k(s,i)}{s}-\pi_k\Big), \ k=m,m+1,...\Bigg)\Rightarrow (Z_k, \ k=m,m+1,...),
\ee

\noindent
where $(Z_k, \ k=m,m+1,...)$ is a mean zero Gaussian process with covariance matrix $V=(v_{rl})_{k\times k}$ given by the limit of the upper left minor size $i\times i$ of the infinite conditional covariance matrix $\hbox{Var}(\d_{D_{s,i}}|\cF_{s,i-1})$, namely
\be{}
\ba{ll}
v_{rr}&=\pi_r(1-\pi_r)=\dfrac{(r+\d)p_r(2m+\d-(r+\d)p_r)}{(2m+\d)^2},\\
v_{rl}&=-\pi_r\pi_l=-\dfrac{(r+\d)(l+\d)}{(2m+\d)^2}p_rp_l, \  r\neq l.
\ea
\ee

\section{Proof of the main result}
\label{proof}

In order to prove Theorem \ref{nodes}, we apply the following martingale central limit theorem.

 \bt{tlc}\cite{RS}
 Let $\{X_{n,m},\cF_{m,n},1\leq m\leq n\}$, $X_{n,m}=(X_{n,m,1},...,X_{n,m,d})^T$ be a $d$-dimensional square-integrable martingale difference array. Consider the $d\times d$ non-negative definite random matrices
 \be{}
 G_{n,m}=(\E[X_{n,m,i}X_{n,m,j}|\cF_{n,m-1}],i,j=1,..,d), \quad V_n=\displaystyle\sum_{m=1}^{n}G_{n,m},
 \ee
 
 \noindent
 and suppose $(A_n)$ is a sequence of $l\times d$ matrices with bounded supremum norm. Assume that
 \bi
 \item[(1)] $A_nV_nA_n^T\stackrel{P}{\ra}\Sigma$ as $n\ra\infty$ for some non-random (hence  non-negative definite) matrix $\Sigma$;
 \item[(2)] $\sum_{m\leq n}\E[X^2\i_{\{{|X_{n,m,i}|>\varepsilon}\}}|\cF_{n,m-1}]\stackrel{P}{\ra}0$ as $n\ra\infty$ for all $i=1,...,d$ and $\varepsilon>0$.
 \ei
 \noindent
 Then, 
 \be{}
 \displaystyle\sum_{m=1}^{n}A_nX_{n,m}\Rightarrow X \quad as \ n\ra\infty,
 \ee
  in $\R^l$, where X is a mean zero l-dimensional Gaussian vector with covariance matrix $\Sigma$.
 \et
We aim at constructing a martingale by appropriate rescaling of the process
\[
N_k(1),N_k(2,1),...,N_k(2,m-1),N_k(2,m),N_k(3,1),...,
\]
for all $k\geq m$. We have the recursion, for $s\geq2$ and $0\leq i\leq m-1$,
\be{ricorsione}
N_k(s,i+1)=N_k(s,i)-\i_{\{D_{s,i+1}=k\}}+\i_{\{k=m\}}\i_{\{i+1=m\}}+\i_{\{k\neq m\}}\i_{\{D_{s,i+1=k-1}\}}.
\ee

The second term on the right-hand side (r.h.s) of (\ref{ricorsione}) accounts for the possibility that the edge being connected at time $i+1$ connects with some vertex of degree $k$. The third term on the r.h.s of \eqref{ricorsione} takes into account the new vertex joining at time $s$, after all of its $m$ edges have been connected, while the fourth term accounts for the possibility that the edge being connected at time $i+1$ might connect with some vertex of degree $k-1$.

If we do not update the degrees during the attachment of a new vertex, but simply at the end of the construction of $\hbox{PA}_s$ starting from $\hbox{PA}_{s-1}$, it is clear that the relation \eqref{ricorsione} changes into
	\be{condizioni}
	\begin{sistema}
		N_k(s,i+1)=N_k(s,i),\ i=0,..,m-2 \\
		\ba{ll}
		N_k(s,m)&=N_k(s-1,m)+\displaystyle\sum_{i=0}^{m-2}\Big[\i_{\{k=m\}}\i_{\{i+1=m\}}\\
		&\quad +\i_{\{k\neq m\}}\i_{\{D_{s,i+1=k-1}\}}-\i_{\{D_{s,i+1}=k\}}\Big].
		\ea
	\end{sistema}
	\ee
For $s\geq2$ and $0\leq i \leq m-1$, let us compute
\be{condexp}
\ba{lll}
\E[N_j(s,i+1)|\cF_{s,i}]&=& N_j(s,i)-\P(D_{s,i+1}=j|\cF_{s,i})+ \i_{\{j=m\}}\i_{\{i+1=m\}}\\
&\quad&+\i_{\{j\neq m\}}\P(D_{s,i+1}=j-1|\cF_{s,i})\\
&=&N_j(s,i)-\frac{j+\d}{s(2m+\d)-2m+i}N_j(s,i)\\
&\quad&+\i_{\{j=m\}}\i_{\{i+1=m\}}+\i_{\{j\neq m\}}\frac{j-1+\d}{s(2m+\d)-2m+i}N_{j-1}(s,i)\\
&=& N_j(s,i)\Big(1-\frac{j+\d}{s(2m+\d)-2m+i}\Big)\\
&\quad&+\i_{\{j\neq m\}}\frac{j-1+\d}{s(2m+\d)-2m+i}N_{j-1}(s,i)+\i_{\{j=m\}}\i_{\{i+1=m\}}.
\ea
\ee
We set
\be{martingala}
M_{s,i}^{(k)}:=a_{s,i}^{(k)}\sum_{j=m}^{k}b_j^{(k)}(N_j(s,i)-\E[N_j(s,i)]), \quad s\geq2, \ 0\leq i \leq m.
\ee
We will show that, for each $k\geq m$, the process $(M_{s,i}^{(k)})$ is a martingale with respect to the same filtration. To this end, we compute
\be{}
\ba{lll}
\E[M_{s,i+1}^{(k)}|\cF_{s,i}]\\
=\E\Bigg[a_{s,i+1}^{(k)}\displaystyle\sum_{j=m}^{k}b_j^{(k)}(N_j(s,i+1)-\E[N_j(s,i+1)])\Bigg|\cF_{s,i}\Bigg]\\
=a_{s,i+1}^{(k)}b_{m}^{(k)}\Biggl(\Bigl(1-\frac{m+\d}{s(2m+\d)-2m+i}\Bigl)N_m(s,i)+\i_{\{i+1=m\}}\\
\quad-\Bigl(1-\frac{m+\d}{s(2m+\d)-2m+i}\Bigl)\E[N_m(s,i)]-\i_{\{i+1=m\}}\Biggl)\\
\quad+a_{s,i+1}^{(k)}\displaystyle\sum_{j=m+1}^{k}b_j^{(k)}\Biggl(N_j(s,i)\Bigl(1-\frac{j+\d}{s(2m+\d)-2m+i}\Bigl)\\
\quad+\frac{j-1+\d}{s(2m+\d)-2m+1}N_{j-1}(s,i)+\Bigl(1-\frac{j+\d}{s(2m+\d)-2m+i}\Bigl)\E[N_j(s,i)]\\
\quad-\frac{j-1+\d}{s(2m+\d)-2m+i}\E[N_{j-1}(s,i)]\Biggl)\\
=a_{s,i+1}^{(k)}\displaystyle\sum_{j=m}^{k}\Biggl(b_j^{(k)}\Bigl(1-\frac{j+\d}{s(2m+\d)-2m+1}\Bigl)\\
\quad+b_{j+1}^{(k)}\frac{j+\d}{s(2m+\d)-2m+1}\Biggl)(N_j(s,i)-\E[N_j(s,i)])\\
\quad+b_k^{(k)}\Bigl(1-\frac{k+\d}{s(2m+\d)-2m+i}\Bigl)(N_k(s,i)-\E[N_k(s,i)]).
\ea
\ee

This is the same as the right-hand side of (\ref{martingala}) if the following holds:
\be{condizioni}
\begin{sistema}
	a_{s,i+1}^{(k)}\Bigl[b_j^{(k)}\Bigl(1-\frac{j+\d}{s(2m+\d)-2m+1}\Bigl)+b_{j+1}^{(k)}\frac{j+\d}{s(2m+\d)-2m+i}\Bigl]=a_{s,i}^{(k)}b_j^{(k)},\ j=m,..,k-1 \\
	a_{s,i+1}^{(k)}b_k^{(k)}\Bigl(1-\frac{k+\d}{s(2m+\d)-2m+i}\Bigl)=a_{s,i}^{(k)}b_k^{(k)}.
\end{sistema}
\ee
Let
\be{coeffa}
\ba{lll}
a_{s,i}^{(k)}&:=&\Biggl[\displaystyle\prod_{t=k-m+1}^{s-1}\prod_{r=0}^{m-1}\Biggl(1-\frac{k+\d}{t(2m+\d)-2m+r}\Biggl)\Biggl]^{-1}\\
&\quad&\times\Biggl[\displaystyle\prod_{r=0}^{i-1}\Biggl(1-\frac{k+\d}{s(2m+\d)-2m+r}\Biggl)\Biggl]^{-1}.
\ea
\ee
Let also
\be{coeffb}
b_j^{(k)}:=\prod_{t=j}^{k-1}\frac{t+\d}{t-k}=(-1)^{k-j}\frac{\G(k+\d)}{(k-j)!\G(j+\d)}, \quad 1\leq j\leq k.
\ee
We adopt the usual convention that $b_k^{(k)}=1$ and $b_j^{(k)}=0$ for all $j>k$.  Direct calculations show that \eqref{coeffa} and \eqref{coeffb} satisfy the conditions in (\ref{condizioni}).

Note that 
\be{ricorrcoeffa}
a_{s,i+1}^{(k)}=a_{s,i}^{(k)}\Biggl(1-\frac{k+\d}{s(2m+\d)-2m+i}\Biggl)^{-1}.
\ee
Moreover, by the Stirling formula, as $s\ra\infty$,
\be{regvar}
\ba{lll}
a_{s,i}^{(k)}&=&\displaystyle\prod_{r=0}^{m-1}\frac{\G\Big(s+\frac{r-2m}{2m+\d}\Big)\Big/\G\Big(k-m+1+\frac{r-2m}{2m+\d}\Big)}{\G\Big(s+\frac{r-k-\d-2m}{2m+\d}\Big)\Big/\G\Big(k-m+1+\frac{r-k-\d-2m}{2m+\d}\Big)}\\
&\quad&\times\displaystyle\prod_{r=0}^{i-1}\Big(1-\frac{k+\d}{s(2m+\d)-2m+r}\Big)^{-1}\\
&\sim& s^{m \frac{k+\d}{2m+\d}}\displaystyle\prod_{r=0}^{m-1}\frac{\G\Big(k-m+1+\frac{r-k-\d-2m}{2m+\d}\Big)}{\G\Big(k-m+1+\frac{r-2m}{2m+\d}\Big)},
\ea
\ee
so that $s\to a_{s,i}^{(k)}$ is regularly varying with index $m(k+\d)/(2m+\d)$.

Thus, we have proved that the process $(M_{s,i}^{(k)})$ is a martingale with respect to the filtration $(\cF_{s,i-1})$. 

\br{compatibilità}
Note that there is not a unique set of coefficients $a_{s,i}^{(k)}$ and $b_j^{(k)}$ that satisfies \eqref{condizioni}. For example, if $a_{s,i}^{(k)},b_j^{(k)}$ is a solution of \eqref{condizioni}, then also $\alpha a_{s,i}^{(k)},b_j^{(k)}$ is a solution for any $\alpha\in\R$. However, if we require that the boundary conditions for $M_{s,i}^{(k)}$ are satisfied, i.e., $M_{s,m}^{(k)}=M_{s+1,0}^{(k)}$ for any $k\geq m$, then \eqref{condizioni} admits a unique solution. This implies that the indexing notations for the processes $N_j(s,i)$ and $M_{s,i}^{(k)}$ are consistent.
\er

There exists a probability mass function $\{p_k, \ k\geq m\}$ such that as $s\ra\infty$ \cite[Chapter 8]{VdH}
\be{conv}
\frac{N_k(s,i)}{s}\ra p_k, \quad k\geq m,
\ee
almost surely, uniformly over $i\in{\{0,...,m-1\}}$.

For $k = m,m+1,...$ and $i\in\{0,...,m-1\}$, define a $k$-variate triangular array of martingale differences by
\be{martdiff}
X_{s,h,j}^{(l)}:=\frac{M_{h,j}^{(l)}-M_{h,j-1}^{(l)}}{a_{s,j}^{(l)}s^{1/2}}, \quad h=k+1,..,s; \ l=m,..,k; \ j=0,..,i.
\ee

In order to use the multivariate martingale central limit theorem, we compute the asymptotic form of the quantities
\begin{align}\label{defg}
	G_{s,h,j}(r,l)&:=\E[X_{s,h,j}^{(r)}X_{s,h,j}^{(l)}|\cF_{h,i-1}]\\
	&=(a_{s,j}^{(r)}a_{s,j}^{(l)}s)^{-1}\E\Bigl[(M_{h,j}^{(r)}-M_{h,j-1}^{(r)})(M_{h,j}^{(l)}-M_{h,j-1}^{(l)})\Bigl|\cF_{h,j-1}\Bigl] \notag
\end{align}
By the martingale property, we get
\be{martprop}
\ba{lll}
\E\Bigl[(M_{h,j}^{(r)}-M_{h,j-1}^{(r)})(M_{h,j}^{(l)}-M_{h,j-1}^{(l)})\Bigl|\cF_{h,j-1}\Bigl] \\

=\E\Biggl[\displaystyle\sum_{d=m}^{r}b_d^{(r)}(a_{h,j}^{(r)}N_{d}(h,j))-a_{h,j-1}^{(r)}N_d(h,j-1))\\
\quad\times\displaystyle\sum_{t=m}^{l}b_t^{(l)}(a_{h,j}N_t(h,j)-a_{h,j-1}N_t(h,j-1))\Biggl|\cF_{h,j-1}\Biggl] \\
\quad-\displaystyle\sum_{d=m}^{r}b_d^{(r)}(a_{h,j}^{(r)}\E[N_d(h,j)]-a_{h,j-1}^{(r)}\E[N_d(h,j-1)])\\
\quad\times\displaystyle\sum_{t=m}^{l}b_t^{(l)}(a_{h,j}^{(l)}\E[N_t(h,j)]-a_{h,j-1}^{(l)}\E[N_t(h,j-1)]).
\ea
\ee
We begin by analyzing the behaviour of the deterministic term on the right-hand side of (\ref{martprop}). When $d=m$, by \eqref{ricorrcoeffa} we obtain
\be{}
\ba{lll}
a_{s,j}^{(l)}\E[N_m(s,j)]-a_{s,j-1}^{(l)}\E[N_m(s,j-1)]\\
=a_{s,j}^{(l)}\Bigl[\i_{\{j=m\}}+\Bigl(1-\frac{m+\d}{s(2m+\d)-2m+j-1}\Bigl)\E[N_m(s,j-1)]\Bigl]-a_{s,j-1}^{(l)}\E[N_m(s,j-1)]\\
=a_{s,j}^{(l)}\Bigl[\i_{\{j=m\}}+\E[N_m(s,j-1)]\frac{l-m}{s(2m+\d)-2m+j-1}\Bigl]\\
\sim a_{s,j-1}^{(l)}\Bigl(\i_{\{j=m\}}+p_m\frac{l-m}{2m+\d}\Bigl),
\ea
\ee
almost surely as $s\ra\infty$, since $a_{s,j}^{(l)}\sim a_{s,j-1}^{(l)}$.

Now consider the case $d\geq m+1$:
\be{}
\ba{lll}
a_{s,j}^{(l)}\E[N_d(s,j)]-a_{s,j-1}^{(l)}\E[N_d(s,j-1)]\\
=a_{s,j}^{(l)}\Bigl(\Bigl(1-\frac{d+\d}{s(2m+\d)-2m+j-1}\Bigl)\E[N_{d}(s,j-1)]+\frac{d-1+\d}{s(2m+\d)-2m+j-1}\E[N_{d-1}(s,j-1)]\Bigl)\\
\quad-a_{s,j-1}^{(l)}\E[N_d(s,j-1)]\\
=a_{s,j}^{(l)}\Bigl(\frac{d-1+\d}{s(2m+\d)-2m+j-1}\E[N_{d-1}(s,j-1)]+\frac{l-d}{s(2m+\d)-2m+j-1}\E[N_d(s,j-1)]\Bigl)\\
\sim \frac{a_{s,j-1}^{(l)}}{2m+\d}\Big((l-d)p_d+(d-1+\d)p_{d-1}\Big),
\ea
\ee
almost surely as $s\ra\infty$. Recall that $p_d$ is given by \cite[(8.42)]{VdH}.

In conclusion, we get
\be{as1}
\ba{lll}
\displaystyle\lim_{s\ra\infty}\frac{1}{a_{s,j-1}^{(l)}}\displaystyle\sum_{d=m}^{l}b_d^{(l)}\Bigl(a_{s,j}^{(l)}\E[N_d(s,j)]-a_{s,j-1}^{(l)}\E[N_d(s,j-1)]\Bigl)\\
=b_m^{(l)}\i_{\{j=m\}}+\frac{1}{2m+\d}\displaystyle\sum_{d=m}^{r}b_d^{(l)}[(l-d)p_d+(d-1+\d)p_{d-1}]
\ea
\ee

Let us focus on the second term of the right-hand side of (\ref{as1}). Rearranging the terms in the sum,
\be{somm0}
\ba{lll}
\displaystyle\sum_{d=m}^{r}b_d^{(r)}[(r-d)p_d+(d-1+\d)p_{d-1}]&=&\displaystyle\sum_{d=m}^{r-1}p_d[(r-d)b_d^{(r)}+(d+\d)b_{d+1}^{(r)}]\\
&=&\displaystyle\sum_{d=m}^{r-1}p_db_d^{(r)}\Bigg[r-d+(d+\d)\frac{b_{d+1}^{(r)}}{b_d^{(r)}}\Bigg]\\
&=&0,
\ea
\ee
since by (\ref{coeffb})
\be{}
 \frac{b_{d+1}^{(r)}}{b_d^{(r)}}=\frac{d-r}{d+\d}.
 \ee
Summarizing, we get
\be{}
\displaystyle\lim_{s\ra\infty}\frac{1}{a_{s,j-1}^{(l)}}\displaystyle\sum_{d=m}^{r}b_d^{(r)}\Bigl(a_{s,j}^{(r)}\E[N_d(s,j)]-a_{s,j-1}^{(r)}\E[N_d(s,j-1)]\Bigl)\\
=b_m^{(r)}\i_{\{j=m\}}.
\ee

Next we look at the first term on the right-hand side of \eqref{martprop}. By (\ref{ricorsione}), we obtain
\be{}
\ba{lll}
a_{s,j}^{(r)}N_d(s,j)-a_{s,j-1}^{(r)}N_d(s,j-1)\\
=a_{s,j}^{(r)}\Bigl(N_d(s,j-1)+\i_{\{d=m\}}\i_{\{j=m\}}+\i_{\{d\neq m\}}\i_{\{D_{s,j}=d-1\}}-\i_{\{D_{s,j}=d\}}\Bigl)\\
\quad-a_{s,j-1}^{(r)}N_d(s,j-1)\\
=a_{s,j}^{(r)}\Bigl(N_d(s,j-1)\frac{r+\d}{s(2m+\d)-2m+j-1}-\i_{\{D_{s,j}=d\}}+\i_{\{d=m\}}\i_{\{j=m\}}\\
\quad+\i_{\{d\neq m\}}\i_{\{D_{s,j}=d-1\}}\Bigl)
:=a_{s,j}^{(r)}\Bigl(N_d(s,j-1)\frac{r+\d}{s(2m+\d)-2m+j-1}+B_s^{(j)}(d)\Bigl),
\ea
\ee
where
\be{defB}
B_s^{(j)}(d):=-\i_{\{D_{s,j}=d\}}+\i_{\{d=m\}}\i_{\{j=m\}}+\i_{\{d\neq m\}}\i_{\{D_{s,j}=d-1\}}.
\ee
Thus, we get
\be{si'}
\ba{lll}
\E\Bigg[\frac{1}{a_{s,j}^{(r)}a_{s,j}^{(l)}}\displaystyle\sum_{d=m}^{r}b_d^{(r)}(a_{s,j}^{(r)}N_d(s,j)-a_{s,j-1}^{(r)}N_d(s,j-1))\\
\quad\times\displaystyle\sum_{t=m}^{l}b_t^{(l)}(a_{s,j}^{(l)}N_t(s,j)-a_{s,j-1}^{(l)}N_t(s,j-1))\Bigg|\cF_{s,j-1}\Bigg]\\
=\E\Bigg[\displaystyle\sum_{d=m}^{r}b_d^{(r)}\Bigg(N_d(s,j-1)\frac{r+\d}{s(2m+\d)-2m+j-1}+B_s^{(j)}(d)\Bigg)\times \\
\quad\displaystyle\sum_{t=m}^{l}b_t^{(l)}\Bigg(N_t(s,j-1)\frac{l+\d}{s(2m+\d)-2m+j-1}+B_s^{(j)}(t)\Bigg)\Bigg|\cF_{s,j-1}\Bigg]
\ea
\ee

By direct calculations, the r.h.s of \eqref{si'} becomes
\be{si}
\ba{lll}
\displaystyle\sum_{d=m}^{r}b_d^{(r)}(r+\d)\frac{N_d(s,j-1)}{s(2m+\d)-2m+j-1}\displaystyle\sum_{t=m}^{l}b_t^{(l)}(l+\d)\frac{N_t(s,j-1)}{s(2m+\d)-2m+j-1}\\
\quad+\displaystyle\sum_{d=m}^{r}b_d^{(r)}(r+\d)\frac{N_d(s,j-1)}{s(2m+\d)-2m+j-1}\displaystyle\sum_{t=m}^{l}b_t^{(l)}\E[B_s^{(j)}(t)|\cF_{s,j-1}]\\
\quad+\displaystyle\sum_{t=m}^{l}b_t^{(l)}(l+\d)\frac{N_t(s,j-1)}{s(2m+\d)-2m+j-1}\displaystyle\sum_{d=m}^{r}b_d^{(r)}\E[B_s^{(j)}(d)|\cF_{s,j-1}]\\
\quad+\displaystyle\sum_{d=m}^{r}\displaystyle\sum_{t=m}^{l}b_d^{(r)}b_t^{(l)}\E[B_s^{(j)}(d)B_s^{(j)}(t)|\cF_{s,j-1}]\\
:=S_{1,s}^{(j)}(r,l)+S_{2,s}^{(j)}(r,l)+S_{3,s}^{(j)}(r,l)+S_{4,s}^{(j)}(r,l),
\ea
\ee
We are interested in the behaviour of the terms $S_{1,s}^{(j)}(r,l)$, $S_{2,s}^{(j)}(r,l)$, $S_{3,s}^{(j)}(r,l)$ and $S_{4,s}^{(j)}(r,l)$ as $s\ra\infty$. We first deal with the terms $S_{1,s}^{(j)}(r,l)$, $S_{2,s}^{(j)}(r,l)$, $S_{3,s}^{(j)}(r,l)$. To this end, note that with probability 1 as $s\ra\infty$
\be{}
\displaystyle\sum_{d=m}^{r}b_d^{(r)}(r+\d)\frac{N_d(s,j-1)}{s(2m+\d)+j-1}\sim\frac{r+\d}{2m+\d}\sum_{d=m}^{r}b_d^{(r)}p_d
\ee
and that
\be{}
\ba{lll}
\displaystyle\sum_{d=m}^{r}b_d^{(r)}\E[B_s^{(j)}(d)|\cF_{s,j-1}]\\
=b_m^{(r)}\E[\i_{\{j=m\}}-\i_{\{D_{s,j}=m\}}|\cF_{s,j-1}]\\
\quad+\displaystyle\sum_{d=m+1}^{r}b_d^{(r)}\E[\i_{\{D_{s,j}=d-1\}}-\i_{\{D_{s,j}=d\}}|\cF_{s,j-1}]\\
=b_m^{(r)}\i_{\{j=m\}}-b_m^{(r)}\P(D_{s,j}=m|\cF_{s,j-1})\\
\quad+\displaystyle\sum_{d=m+1}^{r}b_d^{(r)}(\P(D_{s,j}=d-1|\cF_{s,j-1})-\P(D_{s,j}=d|\cF_{s,j-1}))\\
=b_m^{(r)}\Big(\i_{\{j=m\}}-\frac{m+\d}{s(2m+\d)-2m+j-1}N_m(s,j-1)\Big)\\
\quad+\displaystyle\sum_{d=m+1}^{r}b_d^{(r)}\Big(\frac{d-1+\d}{s(2m+\d)-2m+j-1}N_{d-1}(s,j-1)\\
\quad-\frac{d+\d}{s(2m+\d)-2m+j-1}N_d(s,j-1)\Big)\\
\ea
\ee

Thus, we obtain
\be{}
\ba{lll}
\displaystyle\sum_{d=m}^{r}b_d^{(r)}\E[B_s^{(j)}(d)|\cF_{s,j-1}]\\
\sim b_m^{(r)}\i_{\{j=m\}}+\frac{1}{2m+\d}\displaystyle\sum_{d=m}^{r}b_d^{(r)}((d-1+\d)p_{d-1}-(d+\d)p_d)\\
=b_m^{(r)}\i_{\{j=m\}}-\frac{r+\d}{2m+\d}\displaystyle\sum_{d=m}^{r}b_d^{(r)}p_d,
\ea
\ee
almost surely as $s\ra\infty$, where at the last step we used (\ref{somm0}). Summarizing, with high probability
\be{s123}
\ba{lll}
S_{1,s}^{(j)}(r,l)\ra\frac{(r+\d)(l+\d)}{(2m+\d)^2}\displaystyle\sum_{d=m}^{r}b_d^{(r)}p_d\displaystyle\sum_{t=m}^{l}b_t^{(l)}p_t^,\\
S_{2,s}^{(j)}(r,l)\ra\frac{r+\d}{2m+\d}\displaystyle\sum_{d=m}^{r}b_d^{(r)}p_d\Bigg(b_m^{(l)}\i_{\{j=m\}}-\frac{l+\d}{2m+\d}\displaystyle\sum_{t=m}^{l}b_t^{(l)}p_t\Bigg),\\
S_{3,s}^{(j)}(r,l)\ra\frac{l+\d}{2m+\d}\displaystyle\sum_{t=m}^{l}b_t^{(l)}p_t\Bigg(b_m^{(r)}\i_{\{j=m\}}-\frac{r+\d}{2m+\d}\displaystyle\sum_{d=m}^{r}b_d^{(r)}p_d\Bigg),
\ea
\ee
as $s\ra\infty$. Finally, we consider the term $S_{4,s}(r,l)$. Note that, by (\ref{defB}), we have 
\be{}
\begin{sistema}
B_s^{(j)}(m)=\i_{\{j=m\}}-\i_{\{D_{s,j}=m\}}, \\
B_s^{(j)}(d)=\i_{\{D_{s,j}=d-1\}}-\i_{\{D_{s,j}=d\}}, \quad d\geq m+1.
\end{sistema}
\ee
We separate the analysis of these terms according to the following two events:
\begin{enumerate}
	\item[a)] $\{D_{s,j}=m\}$,
	\item[b)] $\{D_{s,j}=h\}$, $h\geq m+1$.
\end{enumerate}
On the event $\{D_{s,j}=m\}$, by \eqref{defB}, we get
\be{stimab1}
\ba{lll}
\E[B_s^{(j)}(d)B_s^{(j)}(t)|\cF_{s,j-1}]&=&\i_{\{d=t=m\}}\i_{\{j\neq m\}}+\i_{\{d=t=m+1\}}-\i_{\{d=m\}}\\
&\quad&\times\i_{\{t=m+1\}}\i_{\{j\neq m\}}-\i_{\{t=m\}}\i_{\{d=m+1\}}\i_{\{j\neq m\}}.
\ea
\ee
On $\{D_{s,j}=h\}$, $h\geq m+1$ we get
\be{stimab2}
\ba{lll}
\E[B_s^{(j)}(d)B_s^{(j)}(t)|\cF_{s,j-1}]\\
=\i_{\{d=t=m\}}\i_{\{j=m\}}+\i_{\{d=m\}}\i_{\{j=m\}}\Big(\i_{\{t=h+1\}}-\i_{\{t=h\}}\Big)\\
\quad+\i_{\{t=m\}}\i_{\{j=m\}}\Big(\i_{\{d=h+1\}}-\i_{\{d=h\}}\Big)+\i_{\{t=d=h\}}+\i_{\{t=d=h+1\}}\\
\quad-\i_{\{d=h\}}\i_{\{t=h+1\}}-\i_{\{d=h+1\}}\i_{\{t=h\}}.
\ea
\ee

\noindent
By \eqref{si}, combining (\ref{stimab1}) and (\ref{stimab2}) we obtain
\be{s4}
\ba{lll}
S_{4,s}^{(j)}(r,l)&=&\frac{m+\d}{s(2m+\d)-2m+j-1}N_m(s,j-1)\Big[\i_{\{j\neq m\}}(b_m^{(r)}b_m^{(l)}-b_m^{(r)}b_{m+1}^{(l)}\\
&\quad&-b_{m+1}^{(r)}b_m^{(l)})+b_{m+1}^{(r)}b_{m+1}^{(l)}\Big]+\displaystyle\sum_{h=m+1}^{s}\frac{h+\d}{s(2m+\d)-2m+j-1}\\
&\quad&\times N_h(s,j-1)\Big[\i_{\{j=m\}}(b_m^{(r)}b_m^{(l)}+b_m^{(r)}(b_{h+1}^{(l)}-b_h^{(l)})\\
&\quad&+b_m^{(l)}(b_{h+1}^{(r)}-b_h^{(r)}))+b_h^{(r)}b_h^{(l)}+b_{h+1}^{(r)}b_{h+1}^{(l)}-b_h^{(r)}b_{h+1}^{(l)}-b_{h+1}^{(r)}b_h^{(l)}\Big]\\
&\sim&\frac{m+\d}{2m+\d}p_m\Big[\i_{\{j\neq m\}}(b_m^{(r)}b_m^{(l)}-b_m^{(r)}b_{m+1}^{(l)}-b_{m+1}^{(r)}b_m^{(l)})+b_{m+1}^{(r)}b_{m+1}^{(l)}\Big]\\
&\quad&+\displaystyle\sum_{h=m+1}^{\infty}\frac{h+\d}{2m+\d}p_h\Big[\i_{\{j=m\}}(b_m^{(r)}b_m^{(l)}+b_m^{(r)}(b_{h+1}^{(l)}-b_h^{(l)})\\
&\quad&+b_m^{(l)}(b_{h+1}^{(r)}-b_h^{(r)}))+b_h^{(r)}b_h^{(l)}+b_{h+1}^{(r)}b_{h+1}^{(l)}-b_h^{(r)}b_{h+1}^{(l)}-b_{h+1}^{(r)}b_h^{(l)}\Big],
\ea
\ee
almost surely as $s\ra\infty$.

Replacing in \eqref{si} the expressions derived in \eqref{s123} for the terms $S_{1,s}^{(j)}(r,l)$, $S_{2,s}^{(j)}(r,l)$ and $S_{3,s}^{(j)}(r,l)$ and in \eqref{s4} for $S_{4,s}^{(j)}(r,l)$, we conclude that, with high probability,
\be{ajrl}
\ba{lll}
\E\Big[\frac{1}{a_{s,j}^{(r)}a_{s,j}^{(l)}}(M_{s,j}^{(r)}-M_{s,j-1}^{(r)})(M_{s,j}^{(l)}-M_{s,j-1}^{(l)})\Big|\cF_{s,j-1}\Big]\ra a^{(j)}(r,l)\\
:=\frac{m+\d}{2m+\d}p_m\Big(\i_{\{j\neq m\}}(b_m^{(r)}b_m^{(l)}-b_m^{(r)}b_{m+1}^{(l)}-b_{m+1}^{(r)}b_m^{(l)})+b_{m+1}^{(r)}b_{m+1}^{(l)}\Big)\\
\quad+\displaystyle\sum_{h=m+1}^{\infty}\frac{h+\d}{2m+\d}p_h\Big(\i_{\{j=m\}}(b_m^{(r)}b_m^{(l)}+b_m^{(r)}(b_{h+1}^{(l)}-b_h^{(l)})+b_m^{(l)}(b_{h+1}^{(r)}\\
\quad-b_h^{(r)}))+ b_h^{(r)}b_h^{(l)}+b_{h+1}^{(r)}b_{h+1}^{(l)}-b_h^{(r)}b_{h+1}^{(l)}-b_{h+1}^{(r)}b_h^{(l)}\Big)-\Big(b_m^{(r)}\i_{\{j=m\}}\\
\quad-\frac{r+\d}{2m+\d}\displaystyle\sum_{d=m}^{r}b_d^{(r)}p_d\Big)\Big(b_m^{(l)}\i_{\{j=m\}}-\frac{l+\d}{2m+\d}\displaystyle\sum_{t=m}^{l}b_t^{(l)}p_t\Big).
\ea
\ee

Now we consider the quantity

\be{arl}
\ba{lll}
a(r,l)&:=&\displaystyle\sum_{j=0}^{m-1}a^{(j)}(r,l)\\
&=&\frac{m+\d}{2m+\d}p_m\Big((m-1)(b_m^{(r)}b_m^{(l)}-b_m^{(r)}b_{m+1}^{(l)}-b_{m+1}^{(r)}b_m^{(l)})+mb_{m+1}^{(r)}b_{m+1}^{(l)}\Big)\\
&\quad&+\displaystyle\sum_{h=m+1}^{\infty}\frac{h+\d}{2m+\d}p_h\Big(b_m^{(r)}b_m^{(l)}+b_m^{(r)}(b_{h+1}^{(l)}-b_h^{(l)})\\
&\quad&+b_m^{(l)}(b_{h+1}^{(r)}-b_h^{(r)})+m(b_h^{(r)}b_h^{(l)}+b_{h+1}^{(r)}b_{h+1}^{(l)}-b_h^{(r)}b_{h+1}^{(l)}-b_{h+1}^{(r)}b_h^{(l)})\Big)\\
&\quad&-\Big(b_m^{(r)}-\frac{r+\d}{2m+\d}\displaystyle\sum_{d=m}^{r}b_d^{(r)}p_d\Big)\Big(b_m^{(l)}-\frac{l+\d}{2m+\d}\displaystyle\sum_{t=m}^{l}b_t^{(l)}p_t\Big)\\
&\quad&-(m-1)\frac{(r+\d)(l+\d)}{(2m+\d)^2}\displaystyle\sum_{d=m}^{r}b_d^{(r)}p_d\displaystyle\sum_{t=m}^{l}b_t^{(l)}p_t\\
&=&\displaystyle\sum_{h=m}^{\infty}\frac{h+\d}{2m+\d}p_h[(b_m^{(r)}+b_{h+1}^{(r)}-b_h^{(r)})(b_m^{(l)}+b_{h+1}^{(l)}-b_h^{(l)})\\
&\quad&+(m-1)(b_{h+1}^{(l)}-b_h^{(l)})(b_{h+1}^{(r)}-b_h^{(r)})]\\
&\quad&-\Big(b_m^{(r)}-\frac{r+\d}{2m+\d}\displaystyle\sum_{d=m}^{r}b_d^{(r)}p_d\Big)\Big(b_m^{(l)}-\frac{l+\d}{2m+\d}\displaystyle\sum_{t=m}^{l}b_t^{(l)}p_t\Big)\\
&\quad&-(m-1)\frac{(r+\d)(l+\d)}{(2m+\d)^2}\displaystyle\sum_{d=m}^{r}b_d^{(r)}p_d\displaystyle\sum_{t=m}^{l}b_t^{(l)}p_t.
\ea
\ee

By \cite[Section 8]{VdH}, we can write 
\be{defpd}
\ba{lll}
p_d&=&\Bigg((2+\d/m)\frac{\G(m+2+\d+\d/m)}{\G(m+\d)}\Bigg)\frac{\G(d+\d)}{\G(d+3+\d+\d/m)}\\
&:=&c(\d)\frac{\G(d+\d)}{\G(d+3+\d+\d/m)}.
\ea
\ee
Then, by (\ref{coeffb}) and (\ref{defpd}), we get
\be{sommar-1}
\ba{lll}
\displaystyle\sum_{d=m}^{r-1}p_db_d^{(r)}=c(\d)(-1)^r\G(r+\d)\displaystyle\sum_{d=m}^{r-1}(-1)^d\frac{1}{(r-d)!\G(d+3+\d+\d/m)}\\
=c(\d)\frac{(-1)^r\G(r+\d)}{r+2+\d+\d/m}\\
\displaystyle\sum_{d=m}^{r-1}(-1)^d\Bigg[\frac{1}{(r-d)!\G(d+2+\d+\d/m)}+\frac{1}{(r-d-1)!\G(d+3+\d+\d/m)}\Bigg]\\
=c(\d)\frac{(-1)^r\G(r+\d)}{r+2+\d+\d/m}\Big[(-1)^m\frac{1}{(r-m)!\G(m+2+\d+\d/m)}+(-1)^{r-1}\frac{1}{\G(r+2+\d+\d/m)}\Big],
\ea
\ee
where at the last step we used the telescoping property of the sum. Thus by (\ref{sommar-1}), we get
\be{sommar}
\ba{lll}
b_m^{(r)}-\dfrac{r+\d}{2m+\d}\displaystyle\sum_{d=m}^{r}b_d^{(r)}p_d=b_m^{(r)}-\frac{r+\d}{2m+\d}p_r-\frac{r+\d}{2m+\d}\displaystyle\sum_{d=m}^{r-1}b_d^{(r)}p_d\\
=(-1)^{r-m}\frac{\G(r+\d)}{(r-m)!\G(m+\d)}-(-1)^{r-m}\frac{(r+\d)(2+\d/m)\G(r+\d)}{(2m+\d)(r-m)!\G(m+\d)(r+2+\d+\d/m)}\\
=\frac{(-1)^{r-m}\G(r+\d)}{(r-m)!\G(m+\d)}\frac{r+2+\d-r/m}{r+2+\d+\d/m}.
\ea
\ee

\subsection{First condition of Theorem \ref{tlc}}
In this subsection we check that the first condition of Theorem \ref{tlc} holds. We know from (\ref{ajrl}) that
\be{}
sG_{s,s,j}(r,l)\ra a^{(j)}(r,l), 
\ee
almost surely as $s\ra\infty$. From definition (\ref{defg}) we have
\be{gsijh}
G_{s,h,j}(r,l)=\frac{a_{h,j}^{(r)}a_{h,j}^{(l)}}{sa_{s,j}^{(r)}a_{s,j}^{(l)}}hG_{h,h,j}(r,l).
\ee
From the regular variation property (\ref{regvar}), the function
\be{defu}
h\to u(h):=a_{h,j}^{(r)}a_{h,j}^{(l)}hG_{h,h,j}(r,l)
\ee
is regularly varying with index $m(r+l+2\d)/(2m+\d)$. Therefore, from (\ref{gsijh}) and Karamata’s theorem on integration of regularly varying functions \cite[page 25]{R}
\be{kar}
\ba{ll}
\displaystyle\sum_{j=0}^{m-1}V_{s,j}(r,l)&=\displaystyle\sum_{j=0}^{m-1}\displaystyle\sum_{h=k+1}^{s}G_{s,h,j}(r,l)=\displaystyle\sum_{j=0}^{m-1}\displaystyle\sum_{h=k+1}^{s}\dfrac{a_{h,j}^{(r)}a_{h,j}^{(l)}}{sa_{s,j}^{(r)}a_{s,j}^{(l)}}hG_{h,h,j}(r,l) \\
&=\displaystyle\sum_{j=0}^{m-1}\displaystyle\sum_{h=k+1}^{s}\dfrac{u(h)}{sa_{s,j}^{(r)}a_{s,j}^{(l)}} \sim \displaystyle\sum_{j=0}^{m-1}\frac{u(s)}{\Big(m\dfrac{r+l+2\d}{2m+\d}+1\Big)a_{s,j}^{(r)}a_{s,j}^{(l)}} \\ 
&\sim a(r,l)\dfrac{2m+\d}{m(r+l)+2m(\d+1)+\d},
\ea
\ee
for $r,l=m,m+1,...,k$.

This verifies the first condition the martingale central limit theorem of Proposition \ref{tlc} (with each $A_n$ being a $(k-m+1)\times (k-m+1)$ identity matrix).

\subsection{Second condition of Theorem \ref{tlc}}
Next we show that the second condition of Theorem \ref{tlc} holds as well. By (\ref{ricorsione}), we deduce that
\be{}
|(N_l(s,i)-\E[N_l(s,i)])-(N_l(s,i-1)-\E[N_l(s,i-1)])|\leq2 \quad \hbox{for all } l.
\ee

Hence the events $\{|X_{s,i,j,h}^{(l)}|>\varepsilon\}$ are empty for all $l$, for all $h\leq s$ and for all $j\leq i$ as $s\ra\infty$, indeed
\be{}
\ba{lll}
\{|X_{s,i,j,h}^{(l)}|>\varepsilon\}&=&\{|M_{h,j}^{(l)}-M_{h,j-1}^{(l)}|>a_{s,j}^{(l)}s^{1/2}\varepsilon\}\\
&=&\Big\{|a_{h,j}^{(l)}\displaystyle\sum_{t=m}^{l}b_t^{(l)}(N_t(h,j)-\E[N_t(h,j)])-a_{h,j-1}^{(l)}\\
&\quad&\times\displaystyle\sum_{t=m}^{l}b_t^{(l)}(N_t(h,j-1)-\E[N_t(h,j-1)])>a_{s,j}^{(l)}s^{1/2}\varepsilon\Big\}\\
&\subset&\Big\{C_1h^{\frac{l+\d}{2m+\d}}>\varepsilon s^{\frac{l+\d}{2m+\d}}s^{1/2}\Big\}\subset\Big\{C_1>\varepsilon s^{1/2}\Big\},
\ea
\ee
where $C_1>0$ is a large constant independent of $s$. 
Thus, we obtain
\be{}
\i_{\{|X_{s,i,j,h}^{(l)}|>\varepsilon\}}\leq\i_{\{C_1>\varepsilon s^{1/2}\}}\ra0,
\ee
as $s\ra\infty$. This implies that the second condition of Proposition \ref{tlc} holds.

\subsection{The limit process for degree counts}

We conclude that
\be{conv}
\Bigg(\frac{1}{s^{1/2}}\displaystyle\sum_{j=m}^{k}b_j^{(k)}(N_j(s,i)-\E[N_j(s,i)]),\ k\geq m\Bigg)\Rightarrow(Y_k,\ k\geq m)
\ee
in $\mathbb{R}^{\N}$, where $(Y_k, \ k=m,m+1,...)$ is a mean zero Gaussian process with covariance function $R_Y$ given by
\be{}
R_Y(r,l)=\frac{2m+\d}{m(r+l)+2m(\d+1)+\d}a(r,l), \quad r,l\geq m.
\ee
We use this covariance function to define the $(k-m+1)\times(k-m+1)$ matrix with entries
\be{}
R_{Y,k}=(R_Y(r,l), \ m\leq r\leq k, m\leq l\leq k).
\ee
The convergence in (\ref{conv}) means that for all $k=m,m+1,...$
\be{}
C_k\Big(\frac{N_r(s,i)-\E[N_r(s,i)]}{s^{1/2}}, \ r=m,...,k\Big)^T\Rightarrow(Y_r, \ r=m,...,k),
\ee
where $C_k$ is a $(k-m+1)\times(k-m+1)$ matrix with entries
\be{defcrl}
c_{r,l}=
\begin{cases}
	b_l^{(r)} &\hbox{ if } l\leq r, \\
	0 &\hbox{ if } l>r.
\end{cases}
\ee
We are able to conclude that
\be{}
\Big(\frac{N_r(s,i)-\E[N_r(s,i)]}{s^{1/2}}, \ r=m,...,k\Big)\Rightarrow D_k(Y_r, \ r=m,..,k)^T,
\ee
and the covariance matrix of the limiting Gaussian vector is given by
\be{sigk}
\Sigma_k=D_kR_{Y,k}D_k^T.
\ee
Using the following identity, 
\be{ident}
\displaystyle\sum_{t=l}^{r}x^{t-l}b_t^{(r)}b_l^{(t)}=b_l^{(r)}(1+x)^{r-l}, \ 1\leq l\leq r, \ x\in\R,
\ee
it can be shown that $D_k=C_k^{-1}$ has entries
\be{drl}
d_{r,l}=
\begin{cases}
	(-1)^{r-l}b_l^{(r)} &\hbox{ if } l\leq r, \\
	0 &\hbox{ if } l>r.
\end{cases}
\ee

In order to facilitate the computation of the entries of the matrix $\Sigma_k$, by \eqref{arl} we can write $R_{Y,k}$ in the form
\be{ryk}
\ba{lll}
R_{Y,k}&=&\displaystyle\sum_{q=0}^{\infty}h_q\int_{0}^{+\infty}(2m+\d)e^{-[2m(\d+1)+\d]x}R_{q,x}^1dx\\
&\quad&+\displaystyle\sum_{q=m}^{\infty}h_q\int_{0}^{+\infty}(2m+\d)e^{-[2m(\d+1)+\d]x}R_{q,x}^2dx,
\ea
\ee
where 
\be{hq}
h_q=
\begin{cases}
	-1 &\hbox{ if } q=0,1,...m-1, \\
	\frac{q+\d}{2m+\d}p_q &\hbox{ if } q=m,m+1,...,
	\end{cases}
\ee
and the matrix $R_{m,x}$ is the $(k-m+1)\times(k-m+1)$ matrix with entries
\be{rqx}
\ba{ll}
R_{q,x}^1=(C_{q,x}^1)^T(C_{q,x}^1),\\
R_{q,x}^2=(C_{q,x}^2)^T(C_{q,x}^2).
\ea
\ee
The vector $C_{q,x}^1$ has the entries
\be{cqx}
C_{q,x}^1(l)=
\begin{cases}
	\Big(b_m^{(l)}-\frac{l+\d}{2m+\d}\displaystyle\sum_{t=m}^{l}b_t^{(l)}p_t\Big)e^{-mlx}, \  l\geq m &\hbox{ if } q=0, \\
	\Big(\frac{l+\d}{2m+\d}\displaystyle\sum_{t=m}^{l}b_t^{(l)}p_t\Big)e^{-mlx}, \ l\geq m  &\hbox{ if } q=1,...,m-1,\\
	(b_m^{(l)}-b_q^{(l)}+b_{q+1}^{(l)})e^{-mlx} , \ l\geq m &\hbox{ if }  q\geq m.
\end{cases}
\ee
and
\be{cq2}
C_{q,x}^2(l)=
\begin{cases}
0 &\hbox{ if } q=0,1,...,m-1,\\
\sqrt{m-1}(b_{q+1}^{(l)}-b_q^{(l)})e^{-mlx} &\hbox{ if } q\geq m.
\end{cases}
\ee
Therefore
\be{}
\ba{lll}
\Sigma_k&=&\displaystyle\sum_{q=0}^{\infty}h_q\int_{0}^{+\infty}(2m+\d)e^{-[2m(\d+1)+\d]x}(D_k(C_{q,x}^1)^T)C_{q,x}^1D_k^Tdx\\
&\quad&+\displaystyle\sum_{q=m}^{\infty}h_q\int_{0}^{+\infty}(2m+\d)e^{-[2m(\d+1)+\d]x}(D_k(C_{q,x}^2)T)C_{q,x}^2D_k^Tdx.
\ea
\ee

We compute separately the terms in the sum for $q<m$ and $q\geq m$. We begin with $C_{q,x}^1$. For $q\geq m$ and $l=m,m+1,...$, since $b_j^{(k)}=0$ for all $j>k$, by (\ref{ident}) and (\ref{drl}) we have
\be{}
\ba{lll}
(D_k(C_{q,x}^1)^T)(l)\\
=\displaystyle\sum_{t=m}^{l}(-1)^{l-t}b_t^{(l)}(b_m^{(t)}-b_q^{(t)}+b_{q+1}^{(t)})e^{-mtx}\\
=\displaystyle\sum_{t=m}^{l}(-1)^{l-t}b_t^{(l)}b_m^{(t)}e^{-mtx}-\displaystyle\sum_{t=q}^{l}(-1)^{l-t}b_t^{(l)}b_q^{(t)}e^{-mtx}+\displaystyle\sum_{t=q+1}^{l}(-1)^{l-t}b_t^{(l)}b_{q+1}^{(t)}e^{-mtx}\\
=(-1)^{l-m}e^{-m^2x}b_m^{(l)}(1-e^{-mx})^{l-m}-(-1)^{l-q}e^{-mqx}b_q^{(l)}(1-e^{-mx})^{l-q}\\
\quad+(-1)^{l-q-1}e^{-m(q+1)x}b_{q+1}^{(l)}(1-e^{-mx})^{l-q-1}.
\ea
\ee

Therefore, for $q\geq m$ and $r,l=m,m+1,...$, we have

\be{defteta}
\ba{lll}
(D_k(C_{q,x}^1)^T)(l)(D_k(C_{q,x}^1)^T)(r)\\
=(-1)^{l+r}b_m^{(l)}b_m^{(r)}e^{-2m^2x}(1-e^{-mx})^{l+r-2m}\\
\quad-(-1)^{l+r-m-q}(b_m^{(l)}b_q^{(r)}+b_m^{(r)}b_q^{(l)})e^{-m(m+q)x}(1-e^{-mx})^{l+r-m-q}\\
\quad+(-1)^{l+r-q-m-1}(b_m^{(l)}b_{q+1}^{(r)}+b_m^{(r)}b_{q+1}^{(l)})e^{-m(m+q+1)x}(1-e^{-mx})^{l+r-m-q-1}\\
\quad+(-1)^{l+r}b_q^{(l)}b_q^{(r)}e^{-2mqx}(1-e^{-mx})^{r+l-2q}\\
\quad+(-1)^{l+r}(b_q^{(l)}b_{q+1}^{(r)}+b_q^{(r)}b_{q+1}^{(l)})e^{-m(2q+1)x}(1-e^{-x})^{l+r-2q-1}\\
\quad+(-1)^{l+r}b_{q+1}^{(l)}b_{q+1}^{(r)}e^{-2m(q+1)x}(1-e^{-mx})^{l+r-2q-2}\\
:=\displaystyle\sum_{t=1}^{6}\theta_{q,x}^{(t)}(r,l).
\ea
\ee
We have
\be{1}
\ba{lll}
\displaystyle\int_{0}^{+\infty}(2m+\d)e^{-[2m(\d+1)+\d]x}\theta_{q,x}^{(1)}(r,l)dx\\
=(-1)^{l+r}(2m+\d)b_m^{(l)}b_m^{(r)}\displaystyle\int_{0}^{+\infty}e^{-[2m(m+\d+1)+\d]x}(1-e^{-mx})^{l+r-2m}dx\\
=(-1)^{l+r}\dfrac{2m+\d}{m}b_m^{(l)}b_m^{(r)}B(2m+2(\d+1)+\d/m,l+r-2m+1)\\
=(-1)^{l+r}\dfrac{2m+\d}{m}b_m^{(l)}b_m^{(r)}\frac{\G(2(m+\d+1)+\d/m)(l+r-2m)!}{\G(l+r+1+2(\d+1)+\d/m)},
\ea
\ee
where 
\be{beta}
B(\alpha,\beta)=\frac{\G(\alpha)\G(\beta)}{\G(\alpha+\beta)}.
\ee
is the Beta function. Similarly,
\begin{align}
		\notag
	&\displaystyle\int_{0}^{+\infty}(2m+\d)e^{-[2m(\d+1)+\d]x}\theta_{q,x}^{(2)}(r,l)dx\\ \notag
	&=(-1)^{l+r-m-q+1}\dfrac{2m+\d}{m}(b_m^{(l)}b_q^{(r)}+b_m^{(r)}b_q^{(l)})\\
	&\times\frac{\G(m+q+2(\d+1)+\d/m)(l+r-m-q)!}{\G(l+r+1+2(\d+1)+\d/m)}, \\
	 \notag
	&\displaystyle\int_{0}^{+\infty}(2m+\d)e^{-[2m(\d+1)+\d]x}\theta_{q,x}^{(3)}(r,l)dx\\ \notag
	&=(-1)^{l+r-m-q-1}\dfrac{2m+\d}{m}(b_m^{(l)}b_{q+1}^{(r)}+b_m^{(r)}b_{q+1}^{(l)})\\
	&\times\frac{\G(m+q+1+2(\d+1)+\d/m)(l+r-m-q-1)!}{\G(l+r+1+2(\d+1)+\d/m)}, \\
 \notag
&\displaystyle\int_{0}^{+\infty}(2m+\d)e^{-[2m(\d+1)+\d]x}\theta_{q,x}^{(4)}(r,l)dx\\ 
&=(-1)^{r+l}\dfrac{2m+\d}{m}b_q^{(l)}b_q^{(r)}\frac{\G(2q+2(\d+1)+\d/m)(r+l-2q)!}{\G(r+l+1+2(\d+1)+\d/m)}, \\ \notag
&\displaystyle\int_{0}^{+\infty}(2m+\d)e^{-[2m(\d+1)+\d]x}\theta_{q,x}^{(5)}(r,l)dx\\
&=(-1)^{l+r}\dfrac{2m+\d}{m}(b_q^{(l)}b_{q+1}^{(r)}+b_q^{(r)}b_{q+1}^{(l)}\frac{\G(2q+1+2(\d+1)+\d/m)(l+r-2q-1)!}{\G(l+r+2(\d+1)+\d/m)}, \\
 \notag
&\displaystyle\int_{0}^{+\infty}(2m+\d)e^{-[2m(\d+1)+\d]x}\theta_{q,x}^{(6)}(r,l)dx\\
&=(-1)^{l+r}\dfrac{2m+\d}{m}b_{q+1}^{(l)}b_{q+1}^{(r)}\frac{\G(2(q+\d+2)+\d/m)(l+r-2q-2)!}{\G(l+r+2(\d+1)+\d)}.
\end{align}
By (\ref{sommar}) and (\ref{drl}), for $q=0$ and $l=m,m+1,...$ we have
\be{7}
\ba{ll}
(D_k(C_{0,x}^1)^T)(l) \\
=\displaystyle\sum_{t=m}^{l}(-1)^{l-t}b_t^{(l)}\Bigg(b_m^{(t)}-\frac{t+\d}{2m+\d}\displaystyle\sum_{h=m}^{t}b_h^{(t)}p_h\Bigg)e^{-mtx}\\
=\frac{\G(l+\d)}{\G(m+\d)}\displaystyle\sum_{t=m}^{l}(-1)^{t-m}e^{-mtx}\frac{t+2+\d-t/m}{(t-m)!(l-t)!(t+2+\d+\d/m)}.
\ea
\ee

Therefore, by (\ref{7}), for $r,l\geq m$ and $q=0$, we obtain
\be{8}
\ba{lll}
\displaystyle\int_{0}^{+\infty}(2m+\d)e^{-[2m(\d+1)+\d]x}(D_k(C_{0,x}^1)^T)(r)(C_{0,x}^1D_k^T)(l)dx\\
=(2m+\d)\frac{\G(l+\d)\G(r+\d)}{(\G(m+\d))^2}\displaystyle\sum_{t_1=m}^{l}\displaystyle\sum_{t_2=m}^{r}(-1)^{t_1+t_2}\\
\quad\times\frac{[m(t_1+t_2)+2m(\d+1)+\d]^{-1}(t_1+2+\d-t_1/m)(t_2+2+\d-t_2/m)}{(t_1-m)!(t_2-m)!(l-t_1)!(r-t_2)!(t_1+2+\d+\d/m)(t_2+2+\d+\d/m)}
\ea
\ee
Consider now the case $q=1,...,m-1$ and $l=m,m+1,...$. By (\ref{sommar}) and (\ref{drl}), we have
\be{9}
\ba{ll}
(D_k(C_{q,x}^1)^T)(l)&=\displaystyle\sum_{t=m}^{l}(-1)^{l-t}b_t^{(l)}\Bigg(\frac{t+\d}{2m+\d}\displaystyle\sum_{h=m}^{t}b_h^{(t)}p_h\Biggl)e^{-mtx}\\
&=\G(l+\d)\displaystyle\sum_{t=m}^{l}\frac{e^{-mtx}}{(l-t)!\G(t+\d)}\\
&\quad\times\Bigg(b_m^{(t)}-(-1)^{t-m}\frac{\G(t+\d)(t+2+\d-t/m)}{(t-m)!\G(m+\d)(t+2+\d+\d/m)}\Bigg).
\ea
\ee
Therefore, by \eqref{9}, for $r,l\geq m$ and $q=1,...,m-1$, we have
\be{10'}
\ba{ll}
\displaystyle\int_{0}^{+\infty}(2m+\d)e^{-[2m(\d+1)+\d]x}(D_k(C_{q,x}^1)^T)(r)(C_{q,x}^1D_k^T)(l)dx\\
=(2m+\d)\G(l+\d)\G(r+\d)\displaystyle\sum_{t_1=m}^{l}\displaystyle\sum_{t_2=m}^{r}\frac{[m(t_1+t_2)+2m(\d+1)+\d]^{-1}}{(l-t_1)!\G(t_1+\d)(r-t_2)!\G(t_2+\d)}\\
\quad\times\Bigg(b_m^{(t_1)}-(-1)^{t_1-m}\frac{\G(t_1+\d)(t_1+2+\d-t_1/m)}{(t_1-m)!\G(m+\d)(t_1+2+\d+\d/m)}\Bigg)\\
\quad\times\Bigg(b_m^{(t_2)}-(-1)^{t_2-m}\frac{\G(t_2+\d)(t_2+2+\d-t_2/m)}{(t_2-m)!\G(m+\d)(t_2+2+\d+\d/m)}\Bigg).\\
\ea
\ee
By using \eqref{coeffb}, we deduce
\be{10}
\ba{ll}
\displaystyle\int_{0}^{+\infty}(2m+\d)e^{-[2m(\d+1)+\d]x}(D_k(C_{q,x}^1)^T)(r)(C_{q,x}^1D_k^T)(l)dx\\
\quad=\frac{(2m+\d)\G(l+\d)\G(r+\d)}{m^2(\G(m+\d))^2}\displaystyle\sum_{t_1=m}^{l}\displaystyle\sum_{t_2=m}^{r}(-1)^{t_1+t_2}\\
\quad\quad\times\frac{(\d+t_1)(\d+t_2)[m(t_1+t_2)+2m(\d+1)+\d]^{-1}}{(t_1-m)!(t_2-m)!(l-t_1)!(r-t_2)!(t_1+2+\d+\d/m)(t_2+2+\d+\d/m)}.
\ea
\ee
Now we investigate the behavior of the term $C_{q,x}^2$. For $q=m,m+1,...,$, using \eqref{ident} we get
\be{}
\ba{ll}
(D_k(C_{q,x}^2)^T)(l)&=\displaystyle\sum_{t=m}^{l}\sqrt{m-1}(-1)^{l-t}b_t^{(l)}(b_{q+1}^{(t)}-b_q^{(t)})e^{-mtx}\\
&=\sqrt{m-1}\Big((-1)^{q+1-l}e^{-m(q+1)x}b_{q+1}^{(l)}(1-e^{-mx})^{l-q-1}\\
&\quad-(-1)^{q-l}e^{-mqx}b_q^{(l)}(1-e^{-mx})^{l-q}.
\ea
\ee
Thus, we obtain
\be{defteta2}
\ba{ll}
(D_k(C_{q,x}^2)^T)(r)(D_k(C_{q,x}^2)^T)(l)\\
=(m-1)\Big((-1)^{r+l}e^{-2m(q+1)x}b_{q+1}^{(l)}b_{q+1}^{(r)}(1-e^{-mx})^{l+r-2q-2}\\
\quad+(-1)^{r+l}e^{-2mqx}b_q^{(l)}b_q^{(r)}(1-e^{-mx})^{l+r-2q}\\
\quad+(-1)^{r+l}e^{-m(2q+1)x}(b_{q+1}^{(l)}b_q^{(r)}+b_q^{(l)}b_{q+1}^{(r)})(1-e^{-mx})^{l+r-2q-1}\Big)\\
=(m-1)(\theta_{q,x}^{(4)}+\theta_{q,x}^{(5)}+\theta_{q,x}^{(6)}),
\ea
\ee
\noindent
where $\theta_{q,x}^{(i)}$, $i=4,5,6$, are defined in \eqref{defteta}.

Using \eqref{1}-\eqref{8}, \eqref{10} and \eqref{defteta2}, we conclude that the covariance function $R_Z(r,l)$ of the limiting Gaussian process $(Z_k, \ k=m,m+1,...)$ in (\ref{gauss}) is given by (\ref{cov}).

 \end{document}